\documentclass{amsart}[12pt]
\usepackage[]{amsmath, amsthm, amsfonts, verbatim, amssymb}
\usepackage{tikz}
\usetikzlibrary{matrix,arrows}
\usepackage{url}
\usepackage[all,cmtip]{xy}

\usepackage{xcolor}

\usepackage{graphicx}
\usepackage{pinlabel}

\def\oM{\overline M}
\def\cF{\mathcal F}
\def\cH{\mathcal H}
\def\cO{{\mathcal O}}
\def\cI{{\mathcal I}}

\def\cL{{\mathcal L}}

\def\bR{{\mathbb R}}
\def\bZ{{\mathbb Z}}

\def\bC{{\mathbb C}}
\def\bN{{\mathbb N}}

\def\cF{{\mathcal F}}

\begin{document}
\newtheorem {theo}{Theorem}
\newtheorem {coro}{Corollary}
\newtheorem {lemm}{Lemma}
\newtheorem {rem}{Remark}
\newtheorem {defi}{Definition}
\newtheorem {ques}{Question}
\newtheorem {prop}{Proposition}
\def\spb{\smallpagebreak}
\def\mpb{\vskip 0.5truecm}
\def\bpb{\vskip 1truecm}
\def\wtM{\widetilde M}
\def\tM{\widetilde M}
\def\wtN{\widetilde N}
\def\tN{\widetilde N}
\def\tR{\widetilde R}
\def\tC{\widetilde C}
\def\tX{\widetilde X}
\def\tY{\widetilde Y}
\def\tE{\widetilde E}
\def\tH{\widetilde H}
\def\tL{\widetilde L}
\def\tQ{\widetilde Q}
\def\tS{\widetilde S}
\def\tc{\widetilde c}
\def\talpha{\widetilde\alpha}
\def\ti{\widetilde \iota}
\def\hM{\hat M}
\def\hq{\hat q}
\def\hR{\hat R}
\def\bs{\bigskip}
\def\ms{\medskip}
\def\ni{\noindent}
\def\td{\nabla}
\def\pd{\partial}
\def\hol{$\text{hol}\,$}
\def\Log{\mbox{Log}}
\def\bfQ{{\bf Q}}
\def\Todd{\mbox{Todd}}
\def\top{\mbox{top}}
\def\Pic{\mbox{Pic}}
\def\bP{{\mathbb P}}
\def\dxi{d x^i}
\def\dxj{d x^j}
\def\dyi{d y^i}
\def\dyj{d y^j}
\def\dzi{d z^I}
\def\dzj{d z^J}
\def\ozi{d{\overline z}^I}
\def\ozj{d{\overline z}^J}
\def\oz1{d{\overline z}^1}
\def\oz2{d{\overline z}^2}
\def\oz3{d{\overline z}^3}
\def\sI{\sqrt{-1}}
\def\hol{$\text{hol}\,$}
\def\ok{\overline k}
\def\ol{\overline l}
\def\oJ{\overline J}
\def\oT{\overline T}
\def\oS{\overline S}
\def\oV{\overline V}
\def\oW{\overline W}
\def\oY{\overline Y}
\def\oL{\overline L}
\def\oI{\overline I}
\def\oK{\overline K}
\def\oL{\overline L}
\def\oj{\overline j}
\def\oi{\overline i}
\def\ok{\overline k}
\def\oz{\overline z}
\def\om{\overline mu}
\def\on{\overline nu}
\def\oa{\overline \alpha}
\def\ob{\overline \beta}
\def\oGamma{\overline \Gamma}
\def\of{\overline f}
\def\oN{\overline N}
\def\og{\overline \gamma}
\def\ogamma{\overline \gamma}
\def\odelta{\overline \delta}
\def\otheta{\overline \theta}
\def\ophi{\overline \phi}
\def\opd{\overline \partial}
\def\oA{\overline A} 
\def\oB{\overline B}
\def\oC{\overline C}
\def\oD{\overline D}
\def\oIq1{\oI_1\cdots\oI_{q-1}}
\def\oIq2{\oI_1\cdots\oI_{q-2}}
\def\op{\overline \partial}
\def\ua{{\underline {a}}}
\def\us{{\underline {\sigma}}}
\def\Chow{{\mbox{Chow}}}
\def\vol{{\mbox{vol}}}
\def\dim{{\mbox{dim}}}
\def\rank{{\mbox{rank}}}
\def\diag{{\mbox{diag}}}
\def\tor{\mbox{tor}}
\def\supp{\mbox supp}
\def\bp{{\bf p}}
\def\bk{{\bf k}}
\def\a{{\alpha}}
\def\tchi{\widetilde{\chi}}
\def\ta{\widetilde{\alpha}}
\def\ovarphi{\overline \varphi}
\def\ocH{\overline{\cH}}
\def\tV{\widetilde{V}}
\def\tf{\widetilde{f}}
\def\th{\widetilde{h}}
\def\tT{\widetilde T}
\def\hG{\widehat{G}}
\def\hS{\widehat{S}}
\def\hD{\widehat{D}}
\def\Aut{\mbox{Aut}}
\def\hX{\widehat{X}}
\def\Sing{\mbox{Sing}}
\def\tcL{\widetilde{\cL}}
\def\IIm{\mbox{Im}}
\def\alb{\mbox{alb}}
\def\red{\mbox{red}}
\def\Det{\mbox{Det}}
\def\Res{\mbox{Res}}
\def\sing{\mbox{sing}}
\def\tF{\widetilde F}

\ni
\title[Foliations and complex two ball quotients]
{Foliations associated to harmonic maps and some complex two ball quotients}
\author[Sai-Kee Yeung]
{Sai-Kee Yeung}

\address[]{Mathematics Department, Purdue University, West Lafayette, IN  47907USA} \email{yeung@math.purdue.edu
}

\thanks{\noindent{Key words: harmonic maps, foliations, rigidity
The author was partially supported by a grant from the National Science Foundation}}

\ni{\it }
\maketitle
\begin{center}
{\it In memory of Professor Qikeng Lu}
\end{center}

\begin{abstract}  
{\it This article grows out of  some questions about harmonic and holomorphic maps between two bounded symmetric domains.
The main result is an attempt to study foliations associated to a lattice-equivariant harmonic map
of small rank from a complex ball to another.  The result is related to rigidity of some complex ball quotients.  Some open
questions are raised as well.}
\end{abstract}

\bigskip
\begin{center}
{\bf 1. Introduction} 
\end{center} 

\bs\ni {\bf 1.1} 
  In complex geometry, harmonic map is an important tool in the study of rigidity in the last few decades.  
A celebrated result of  
Siu [Siu1] [Siu2] states that a harmonic map $\Phi:M\rightarrow N$ between K\"ahler manifolds with $\rank_{\bR}(\Phi)$
 sufficiently large is holomorphic or conjugate-holomorphic as a result of $\pd\op$
Bochner formula.  A natural question is to study the same problem when $\rank_{\bR}(\Phi)$ is smaller than a critical value.

The problem is particularly interesting in the special case that 
$N$ is a Hermitian symmetric space of non-compact type.  In such 
case, the universal covering $\tN$ of $N$ is also a bounded symmetric domain and is classified by E. Cartan.  We know that there are six types of bounded
symmetric domains denoted by $D^I_{m,n}, D^{II}_n, D^{III}_n, D^{IV}_n, D^V$ and $D^{VI}$.  In Theorem 6.7 of [Siu2], it is
proved that a harmonic map $\Phi:M\rightarrow N$ with $\rank_{\bR}(\Phi)\geqslant 2p(N)+1$ is either holomorphic or
anti-holomorphic, where
$p(D^I_{m,n})=(m-1)(n-1)+1,$ $p(D^{II}_{n})=\frac12(n-2)(n-3)+1,$ $p(D^{III}_{n})=\frac12 n(n-1)+1,$ $p(D^{IV}_{n})=2,$
$p(D^{V})=6,$ and $p(D^{IV})=11.$

Hence there is the following natural question.
\begin{ques}
Let $\Phi:M\rightarrow N$ be a harmonic map from a K\"ahler manifold to a Hermitian symmetric space $N$.
Suppose that $\rank_{\bR}(\Phi)<2p(N)+1.$  Under what conditions can we conclude that $\Phi$ is either holomorphic 
or anti-holomorphic?
\end{ques}

 The answer is not clear even
if both of the domain and image of the manifolds are complex balls or complex ball quotients.  Complex ball of complex dimension $B^n_{\bC}$
are identified as $D^I_{1,n}$.  Hence for $n=2$, $2p(B^2_{\bC})+1=3$ and we only need to worry about the cases that
$\rank_{\bR}(\Phi)\leqslant 2$ from the results of Siu mentioned above.
 In such case,  there are results which show
in various situations that
a fibration with compact leaves exists $M$, such as in  [Siu3] and [CT].  We refer the reader to [Mo1] and
[JY] for other results related to fibrations coming from harmonic maps.  

Suppose we consider now that $M$ is a compact manifold or a non-compact manifold of finite volume.  Let  $\tM$ be the
universal covering of $M$.
Consider now a $\pi_1(M)$ equivariant harmonic map from $\tM$ to $\tN$.  Then even if we have
an answer for Question 1, a holomorphic fiber on $\tM$ may not descend to a subvariety on $M$, but rather a foliation,
possibly with singularities.
A specific question is the following.

\begin{ques}
Let $\tM$ and $\tN$ be complex
balls possibly of different dimension.
Let $\Gamma$ be a cocompact lattice in $\Aut(\tM)$, the automorphism group of $\tM$.  Let $\rho:\Gamma\rightarrow \Aut(\tN)$ be
a homomorphism.
Let $\Phi:\tM\rightarrow \tN$ be a $\rho$-equivariant harmonic map with small $\rank_{\bR}(\Phi)$. Does a generic fiber of  $\Phi$ descends to 
complex subvariety or foliation on $M$?
\end{ques}

\ni {\bf 1.2} In this paper we only investigate the simplest such situation in which the domain and the image $M$ and $N$ are complex balls with $\dim_{\bR}(\tM)=2$ and  $\dim_{\bR}(\Phi(\tM))=1$ or $2$ respectively.  In fact, we would study in details only the case that $M$ has the smallest 
possible Euler number, namely $3$.

\begin{theo}
Let $M$ be a compact complex two ball quotient with $c_2(M)=3$. 
 Let $\Phi:\tM\rightarrow B_{\bC}^1$ be a non-constant $\Gamma$-equivariant holomorphic map.
Then the fibers of $\Phi$ descends to compact algebraic curves on $M$.  Moreover, the fibers comes from 
the fibers of a fibration $f:M\rightarrow R$ to a projective algebraic curve $R$.
\end{theo}

The proof actually depends on two steps.  First of all, using results of Simpson [Sim] on complex variation of Hodge structure, see also related earlier work of [Siu3], [CT] and [JY], we
prove that the fiber of $\Phi$ on $\tM$ are analytic curves.  The fibration on $\tM$ descends to gives rise to a holomorphic foliation on $M$.
The second step is the main task in this paper.  We need to show that a leaf of the foliation is compact to arrive a contradiction.  This requires
a delicate analysis of the associated holomorphic foliation, possibly with singularities, on a complex two ball quotient.
As an application, we obtain the following result.

\begin{theo}
Let $M=B_{\bC}^2/\Gamma$ be a compact complex two ball quotient with $c_2(M)=3$.  Denote by $\tM\cong B_{\bC}^2$ the universal covering of
$M$.  Let $\rho:\Gamma\rightarrow PU(2,1)$ be a rigid reductive representation.  Then there exists a $\Gamma$-equivariant holomorphic map
$\Phi:\tM\rightarrow B_{\bC}^2$
of complex rank $2$. 
\end{theo}

Theorem 2 has applications to problems in rigidity for complex hyperbolic lattices of small covolume.  In fact, it provides some missing details in the proof of Archimedean rigidity in [Y3].
We refer the readers to \S4.2 for more details.

\ms
\ni {\bf 1.3} Here is the organization of the paper.  In \S2, we recall some basic facts about complex two ball quotients of small Euler numbers.
\S3 is the main part of this paper, where we study a foliation on $M$ derived from a harmonic map, and prove that actually leaves
of the foliation are compact.  The result was used in \S4 to give some conclusions about rigidity properties.  \S4 also gives
correction of some statements in [Y3] and provides some details for the arguments in [Y3].  Some concluding remarks and problems related to questions raised in this article are given in  \S5.

\ms
It is a pleasure for the author to thank Domingo Toledo for helpful discussions and comments.

\bs
\begin{center}
{\bf 2. Preliminaries on ball quotients of small Euler numbers} 
\end{center} 

\bs\ni
{\bf 2.1}  Since we are mainly concerned with complex two ball quotients with small Euler number in this article, we collect  in this section
some general facts about the Euler number of complex two ball quotients.   Most of the results mentioned in this section are
known.  By a complex two ball quotient, we are considering $M=B_{\bC}^2/\Gamma=\Gamma\backslash PU(2,1)/P(U(2)\times U(1)$,
where $\Gamma$ is a torsion-free lattice of $PU(2,1)$, which is the automorphism group of $B_{\bC}^2$.  In particular, $M$ is a smooth
projective algebraic manifold.  From explicit Chern form computations, there is the well-known Chern number identity
$c_1^2=3c_2$, where the second Chern number $c_2$ is the same as the Euler number of $M$.

\begin{prop} Let $M$ be a compact complex two ball quotient.\\
(a) The smallest Euler number $c_2(M)$ of a compact complex two ball quotient is $3$.\\
(b) There exists a sequence of complex ball quotients $M_i$ with $c_2(M_n)=3n$ for all $n\in \bN$, the set of natural numbers.\\
(c) We may find $M_n$ in (b) such that  $b_1(M_n)>0$ for each $n>0.$ 
\end{prop}

\ni{\bf Proof} (a) can be deduced easily from Riemann-Roch Formula.  We refer the readers to [Y2] or [Y3] for some elaborations.

For (b), there exist smooth complex ball quotients $M$ with $c_2(M)=3$ provided by fake projective planes, which are compact
complex two ball quotients with Betti numbers the same as the projective plane and hence the first Betti number is trivial.  Cartwright and Steger
show in [CS] that there exists a compact complex ball quotients with $c_2(M)=3$ and $h^1(M)=1$.  Since $H_1(M,\bZ)$ is the abelianization 
of $\pi_1(M)$, we have the following
 homomorphism.
\begin{equation}
f_n:\Gamma\rightarrow H_1(M)=\bZ\rightarrow \bZ/n\bZ. \label{eq:ab-cover}
\end{equation}
Now $\Gamma_n:=\ker(f_n)$ is a lattice from which $M_n:=B_{\bC}^2/\Gamma_n$ satisfies $c_2(M_n)=3n$.  The above argument is given in [CS].

Consider now (c).  Since $b_1(M)=2$ for the Cartwright-Steger surface $M$, holomorphic one forms on $M$ can be pulled back 
by the natural projection map to give non-trivial holomorphic
one forms on $M_n$ for each $M_n$ mentioned in (b).

\qed

\bs\ni
{\bf 2.2}  We remark here that there is a corresponding statement for quotients of complex two ball $B_{\bC}^2$ by cofinite neat lattices
in $PU(2,1)$.   We do not consider cofinite lattices in this article, but clearly Questions 1 and 2 can be phrased for such
locally Hermitian symmetric spaces as well.  

Recall that a lattice $\Gamma$ is said to be neat if the subgroup in $\bC^*$ generated by eigenvalues of elements
of $\Gamma$ as matrices in general linear group is torsion-free.
In such case, we define the Chern numbers $c_1^2$ and $c_2$ as the integral of the standard Chern forms in terms of 
the Poincar\'e metric on $M$ over the manifold .

\begin{prop} Let $\Gamma$ be a cofinite neat lattice of $PU(2,1)$.  Let $M=\tM/\Gamma$ denote
the corresponding smooth complex two ball quotient of finite volume.\\
(a) The smallest possible Euler number $c_2(M)$ is $1$. \\
(b) $c_2(M)=1$ implies that $h^1(M)>0$. \\
(c) There exists a sequence of complex ball quotients $M_n$ with $c_2(M_n)=n$ for all $n\in \bN$, the set of natural numbers.\\
(d) We may find $M_n$ in (c) such that  $b_1(M_n)>0$ for each $n>0.$ 
\end{prop}

\ni{\bf Proof} (a)  Since $M$ is a complex two ball quotient, by looking at the Chern forms, we know that
$c_2=\int_M C_2=\frac13\int_M C_1\wedge C_1>0.$  It is well-known from [AMRT] and [Mo2] that there exists a
 smooth toroidal compactification $\oM=M\cap_i T_i$ of $M$ obtained by attaching torus or tori $T_i$ to the ends of $M$,where $T_i$ are
 disjoint.  It is also well-known
 that $c_2(M)=c_2(\oM)-\sum_i e(T_i)=c_2(\oM)$ are integers,  where $e(T_i)$ denotes the Euler number of $T_i$.  Hence $c_2(M)$ is a positive
 integer and is $\geqslant 1$.  $T_i\cdot T_i<0$ as they are contractible.
 A cofinite neat arithmetic lattice $\Gamma$ with $c_2(M)=1$ was given by Parker [P], see also [H].

For part (b), we denote $D=\sum_i T_i$ the compactification divisor.   Let $K=K_{\oM}$ be the canonical line bundle of $\oM$.  
 It is well known from the growth of the Poincar\'e metric near $D$, cf. [Mu],  that 
 $$c_1(M)^2=c_1(\oM,K+D)^2=K^2+2K\cdot T_i+\sum_i T_i^2=K^2-K\cdot T_i,$$
 where in the last step we used the fact that $K\cdot T_i+T_i\cdot T_i=e(T_i)=0$ as $T_i$ is an elliptic curve, where $e(T_i)$ is the
 Euler number of $T_i$.
 
We conclude from Noether's Formula that
$$h^0(\oM)-h^1(\oM)+h^2(\oM)=\frac1{12}(K^2+c_2(\oM))=\frac1{12}(c_1^2+\sum_iT_i^2+c_2(M))=\frac{c_2}3+\frac{\sum_i T_i^2}{12}.$$
Since $c_2(M)=1$, $T_i\cdot T_i<0$ and $h^0(\oM)=1$, we conclude that $h^1>0$.

(c) follows by applying (\ref{eq:ab-cover}) to an example of $c_1(M)=1$ in (a) and (b) with $M_1=M$, namely, letting
$\Gamma_n=\ker(f_n)$ for $f_n:\Gamma\rightarrow H_1(M)\rightarrow \bZ\rightarrow \bZ/n\bZ.$

  (d) follows by pulling back holomorphic one forms
from $M_1=M$.

\qed

\begin{center}
{\bf 3.  Foliation and harmonic map} 
\end{center} 

\bs\ni
{\bf 3.1}  Throughout the section, we assume that $M=B_{\bC}^2/\Gamma$ is a compact complex two ball quotient.  $\tM=B_{\bC}^2$ is the universal
covering of $M$.  Let $\Phi:\tM\rightarrow B_{\bC}^1$ be a non-constant holomorphic map.  
Observe now that if a fiber $\eta$ of $\Phi$ has multiplicity more than $1$, the image of the fiber $\eta$ must be a closed curve in $M$. In fact, if the image of $\eta$ is not a closed curve on $M$, it would have some limit points in $M$ along some local transverse cross section to the foliation. This implies that 
$\Phi$ has multiplicity greater than $1$ for a generic fiber of $\Phi$, a contradiction. Hence multiple fibers correspond to compact invariant curves of the vector field on $M$, and there are at most a finite number of those leaves, denoted by $E_1,\dots,E_k$.

We observe that $\ker(d\Phi)$ defines a foliation $\cF$ on $\tM$ and hence on $M$, defined by $\Phi^*dw=0$, where $w$ is a local holomorphic coordinates at a point on $\Phi(\tM)$.  
The expression $d\Phi$ is non-degenerate on $M \backslash \cup_{i=1}^k E_i $ apart from a finite number of points. The foliation on  $M \backslash \cup_{i=1}^k E_i $ is non-degenerate generically and extends naturally across $\cup_{i=1}^k E_i$ to give a foliation $\cF$ with a finite number of singularities on $M$. This is essentially the same as the foliation obtained from $d\Phi$, but neglecting the multiplicities of the closed curves $\cup_{i=1}^k E_i$. In other words, we consider the saturation of a foliation if it is defined locally by holomorphic one forms or vector fields as described in page 11 of [Br2]. Hence we denote by $\Sing(\cF)$, the singularity set of $\cF$, a discrete set of points.

We refer for example to [Br1] for standard notations about holomorphic foliations on a complex surface.
There is a short
exact squence 
\begin{equation}
0\rightarrow T_\cF\rightarrow T_M\rightarrow \cI_ZN_\cF\rightarrow 0
\end{equation}
which is dual to 
\begin{equation}
0\rightarrow N^*_\cF\rightarrow \Omega_M\rightarrow \cI_ZT^*_\cF\rightarrow 0
\end{equation}
after tensoring with $K_M$, where $N_{\cF}$ is the normal bundle to the foliation, and $\cI_Z$ is the ideal sheaf with support on the singularity of the foliation.

There is an associated long exact sequence
$$\begin{array}{ccccccccc}
0&\rightarrow&H^0(M,N^*_{\cF})&\stackrel{\iota}\rightarrow&H^0(M,\Omega_M)&\stackrel{\alpha}\rightarrow&H^0(M,\cI_ZT^*_{\cF})&&\\
&\stackrel{\beta}\rightarrow&H^1(M,N^*_{\cF})&\stackrel{\gamma}\rightarrow&H^1(M,\Omega_M)&\stackrel{\delta}\rightarrow&H^1(M,\cI_ZT^*_{\cF})&&\\
&\stackrel{\epsilon}\rightarrow&H^2(M,N^*_{\cF})&\stackrel{\sigma}\rightarrow&H^0(M,\Omega_M)&\stackrel{\tau}\rightarrow&H^0(M,\cI_ZT^*_{\cF})&\rightarrow&0.
\end{array}$$

\ms\ni
{\bf 3.2}  
\begin{lemm} Let $\Phi:\tM\rightarrow B_{\bC}^1$ be a non-constant $\Gamma$-equivariant holomorphic map as above.  Then unless $h^0(M,N^*_{\cF})=0$ and $h^0(M,\Omega_M)\leqslant 1$, the fibers of $\Phi$ descend to hyperbolic projective algebraic curves on $M$.
\end{lemm}

\ni{\bf Proof}  \ For $h^0(M,N^*_{\cF})\geqslant 2$, the two linearly independent elements in $\iota(H^0(M,N^*_{\cF})\subset H^0(M,\Omega_M)$ and the Castelnouvo-de Franchi argument, cf. [BHPV],
lead to the conclusion that the foliation comes from a holomorphic fibration over a curve of genus at least $2$.  In such case, the fibers are all curves on $M$ and 
we are done.

 Suppose $h^0(M,N^*_{\cF})=1$. Let $\omega\in  H^0(M,N^*_{\cF})$.    Denote by the same symbol
its pull-back to $M$. From definition, $\omega$ annihilates the tangent vectors to fibers of
$\Phi$ on $M$. We claim that any leaf of $\cF$ must be compact. Suppose on the contrary
that a leaf $\cL$ of $F$ is not compact on $M$. Denote by $\tcL$ a lift of $\cL$ in $\tM$ . $\tcL$ is the   
fiber of $\Phi:M\rightarrow C$.   Suppose $\Phi(\tcL)=o\in C$.  Let $V$ be a small coordinate neighborhood of $o$  in $C$ and $z$ be a coordinate function on $V$ . Then 
 $\Phi^{-1}(V)$ is a neighborhood of $\tcL$. As both $\pi^*\omega$ and $\Phi^*dz$ annihilate tangent vectors to $\tcL$ on $M$ , we know that in a neighborhood of $\tcL$, 
 $\pi^*\omega= f\Phi^*dz$ for some function $f$ which is holomorphic along $\tcL$. Taking exterior derivative
 $$0=d\pi^*\omega=df\wedge\Phi^*dz+fd\Phi^*dz.$$
It follows that $df(v) = 0$ for all tangent vectors $v\in T\tcL$. Hence $f$ is constant along $\tL$. In other words, we may regard $\pi^*\omega=\Phi^*dz$ along $\tcL$. Note that the same argument shows that $\pi^*\omega|_{\tcL}=\Phi^*\eta|_{\tcL}$ for any local holomorphic one form $\eta$ on $V$ as long as $\Phi^*\eta|_{\tcL}=\Phi^*dz|_{\tcL}$.  In particular, this holds for $\eta= h(z)dz$, where $h$ is a holomorphic function on $V$ with $h(o) = 1$. 

Denote $U=\Phi^{-1}(V)$.  Let $W$ be a small coordinate neighborhood of $\tM$ which has non-empty intersection with $U$. Since $\cL$ is dense on $M$, we may assume that the image $\pi(U\cap W)$ contains an infinite number of disconnected
pieces $\cL_i$, $i\in \bN$ of $\cL\cap\pi(U\cap W)$, by taking a smaller $W$ if necessary. 
Denote $\tL_i=\pi^{-1}(\tL_i)$ on $\Phi^{-1}(V)\cap W$.
From the above discussions, $\omega|_{\tcL_i}=\Phi^*dz|_{\tcL_i}$ for each $i$ on $U\cap W\subset \tM$, where
we identify $\pi^{-1}(\cL_1)$  with $\tcL$ on $W$.  As remarked above, the same argument shows that $\pi^*\omega|_{\tcL_i}=\Phi^*(h(z)dz)|_{\tcL_i}$.  Since we may choose $h(z)$ so that $h(z)\neq 1$ on $\tcL_2$, we immediately reach a contradiction. Hence all  the leaves of $\cF$ on $M$ are compact. It follows that $\Phi$ induces a fibration $p : M\rightarrow R$ to an algebraic curve $R$ of genus at least $1$, as in the proof of the classical Castelnouvo-de Franchi Theorem mentioned above, see also Proposition 6.2 of [Br2]. The fibration lifts to a fibration $\widehat{p}:\tM\rightarrow R$.  Since the fibration in this case is actually the original one induced from $\Phi$, which is $\Gamma$-equivariant, it follows that there is an induced covering map $q:C=\Phi(\tM)\rightarrow R$.  As $\Gamma$ acts freely on $M$, $\Phi_*\Gamma$ acts faithfully on $C$ as well. We conclude that the Poincar\'e metric of $C$ descends to $R$. Hence $R$ is hyperbolic and has to be of genus at least $2$.  Again we are done in this case.

Consider now the case that $h^0(M,N^*_{\cF})=0$.  Assume that $h^0(M,\Omega_M)\geqslant 2$. From the long exact sequence 
studied above, we know that $h^0(M,\cI_ZT^*_{\cF})-\dim\IIm(\beta)\geqslant 2$.  Hence we may find two linearly independent elements $\omega_1, \omega_2\in H^0(M,\Omega_M)\cong \IIm(\alpha)\subset H^0(M,\cI_ZT^*_{\cF})$.  As $T^*_{\cF}$ is one dimensional, we may write $\omega_1=f\omega_2$.  Clearly $\omega_1\wedge\omega_2=0$.  Castelnouvo-de Franchi Theorem implies that there is a fibration of $M$ over a curve $R$ of genus at least $2$.  This concludes the proof of Lemma 1.

\qed

\ms\ni
{\bf 3.3}  We may now proceed to the proof of our main theorem.

\ms
\ni{\bf Proof of Theorem 1}  Since $M$ is a compact complex two ball quotient, we know that $c_1^2=3c_2=9$ and arithmetic genus
$\chi(\cO)=\frac12(c_1^1+c_2)=1$.
  From the above lemma, it suffices for us to assume that $h^0(M,\Omega_M)\leqslant 1$.  
  
  Suppose
$h^0(M,\Omega_M)=0$.  It follows from Hodge Theory that the first Betti number $b_1=0$ and hence the third Betti number $b_3=0$.
From $h^{0}(\cO)-h^1(\cO)+h^2(\cO)=\chi(\cO)=1$ we conclude that $h^2(\cO)=0$.  
From Riemann-Roch, $c_2=\sum_{i=0}^4(-1)^ib_i$ and hence the second Betti number $b_2=1$.  It follows from Hodge Decomposition that
$h^{1,1}=1$.  
Let $\omega_C$ be the K\"ahler form of the Poincar\'e metric on $C$.  Since the pull-back $\Phi^*\omega_{C}$ is $\Gamma$-invariant
and descends to $M$, it has to be a non-trivial multiple of $\omega_M$, the standard K\"ahler form on $M$.  Here we use the fact that $h^{1,1}=1$.
This however contradicts the fact that 
$\pi^*\omega_C\wedge\pi^*\omega_C=0$ as $C$ has dimension $1$.

\ms
Hence we may assume that $h^0(M,\Omega_M)=1$ and $h^0(M,N^*_{\cF})=0$.  Similar to the discussions in the last paragraph, if $h^{1,0}=1$, we would have $h^{1,1}=3$.
From the long exact sequence in {\bf 3.1}, it follows that 
$h^0(M,\cI_ZT^*_{\cF})-\dim(\IIm(\beta))=1$.  In this case, there is a non-trivial Albanese map $\alb:M\rightarrow E$,  an elliptic curve generated by the holomorphic one form $\theta\in H^0(M,\Omega_M)$.  As $\theta$ can be considered as an element in $H^0(M,\cI_ZT^*_{\cF})$, the fibers of $\alb$ is
actually transversal to the leaves of $\cF$ generically.  Recall that we have the following Euler number identity for fibration as given in [BPHV], pp 157,
\begin{equation}
e(M)=e(E)e(M_s)+\sum_{s_o}n_{s_o}, \label{eq:euler}
\end{equation}
where $M_s$ is a generic fiber and the sum is taken over the finite number of singular fibers $M_{s_o}$ of $\pi$, of 
which each $n_{s_o}=e(M_{s_o})-e(M_s)$ is a non-negative integer, and is positive unless $M_{s_o}$ is a multiple fiber with
$(M_{s_o})_{\red}$ a nonsingular elliptic curve (cf. [BHPV], page 118).

Applying the above identity to our Albanese fibration, we know that the contribution of the singular set
$Z_{\alb}$ of $\alb$ is given by $\sum_{i\in Z_{\alb}} \delta_{\alb,i}=3$.  Since $\theta$ can be considered to be living in $H^0(M,\cI_ZT^*_{\cF})$,
we know that the singularity set of the foliation $Z\subset Z_{\alb}$.

We now consider the singularities of $\cF$. Recall an invariant of foliation introduced by Baum-Bott in [BB], see also [Br1] for details in the following setting. In a neighborhood of a singular point $p\in \Sing(\cF)$, we may suppose that $\cF$ is generated
by a vector field $v$ given in local coordinates by $v(z,w)=F(z,w)\frac{\pd}{\pd z}+G(z,w)\frac{\pd}{\pd w}$.  Then
$$\Det(p,\cF)=\Res_0\{\frac{\det J(z,w)}{F(z,w)G(z,w)} dz\wedge dw\},$$
is a non-negative integer. Define
$$\Det(\cF)=\sum_{p\in\sing(\cF)}\Det(p,\cF).$$
From [BB], see also Proposition 1 of [Br1], we know that
\begin{equation}
\Det(\cF)=c_2(M)-c_1(T_{\cF})\cdot c_1(M)+c_1^2(T_{\cF}). \label{eq:a}
\end{equation}
It follows from the equation that 
\begin{equation}
c_1(T_{\cF})\cdot c_1(N_{\cF})=c_2(M)-\Det(\cF). \label{eq:b}
\end{equation}

Note that apart from the finite number of singular points, the normal bundle $N_{\cF}$ can be locally represented by the lift of the tangent vectors to the base of the fibration $\Phi:\tM\rightarrow C$.  From Riemann-Roch formula for $T_{\cF}$ on M, it follows that
$c_1(T_{\cF})\cdot c_1(N_{\cF})=c_1(T_{\cF})\cdot (c_1(T_{\cF})-c_1(T_M))$ is an even integer. Hence as 
$c_2(M)=3$ and $\Det(\cF)\geqslant 0$, identity (\ref{eq:b}) implies that $c_1(T_{\cF})\cdot c_1(N_{\cF})=2$ and $\Det(\cF)=1$, the latter implies that there is only one singularity $Q$ for $\cF$, with Milnor number $1$.

From the earlier discussions, we conclude that $Q\in Z_{\alb}$.  As $\sum_{i\in Z_{\alb}} \delta_{\alb,i}=3$ and the singularity of $Q$ can only contribute $1$
to the sum, we conclude that there is at least one more point, say $Q'\in Z_{\alb}\backslash Z_{\cF}$.  Since by construction, the foliation is smooth around $Q$, we may choose a local  
coordinate system $(x,y)\in W$ centered at $Q'$ so  that leaves of the foliation are given by $y=c$, a small constant.
In such a coordinate, we may write $\theta= fdx$ for some local holomorphic function $f$ on $W$.   As $\theta$ vanishes at $Q'$,   it follows that $f$ has a non-trivial zero divisor passing through $Q'$ on $W$.  Since $\theta$ is is the pull-back of a one form on the elliptic curve E by the Albanese map 
$\alb:M\rightarrow E$, this is possible only if alb contains a singular fiber. In other words, $\theta$ contains a multiple fiber. 
 
 Recall now the expression (\ref{eq:euler}) with respect to the Albanese fibration. Consider first the contributions from multiple fibers. For a multiple fiber $s_o$,
the expression $n_{s_o}$ is positive  unless the reduced $M_{s_o}$ is an elliptic curve, which is not possible as $M$ is hyperbolic. In such a case, the fact that we have a small Euler number $3$ implies that a multiple fiber $M_{s_o}$ has just one reduced irreducible component $D_{s_o}$ and we may write
$M_{s_o}=\alpha_{s_o}D_{s_o}$ for some $\alpha_{s_o}\geqslant 2$, cf. Corollary 2 in {\bf 5.3} of [CKY].  Note that $n_{s_o}=2(g(M_{s_o})-g(M_s))$ and 
$2(g(M_s)-1)=2\alpha_{s_o}(g(D_{s_o})-1)$.  It follows that $n_{s_o}=2(\alpha_{s_o}-1)(g(M_{s_o})-1)$.  Since $g(D_{s_o})\geqslant 2$ and $\alpha_{s_o}\geqslant 2$,
we conclude from (\ref{eq:euler}) that there is precisely one multiple fiber $M_{s_o}$, and furthermore, $g(D_{s_o})=2$ and $\alpha_{s_o}=2$.  There is also the contribution of $Q$ to the formula (\ref{eq:euler}). The point $Q$ cannot lie on $M_{s_o}$, for otherwise the normalization of
$D_{s_o}$ would have genus at most $1$ and would not be hyperbolic.

It follows from our setting that $\theta/\xi$ is precisely an element in $H^0(M,\cI_ZT^*_{\cF})$, where $\xi$ is the canonical section of $D_{s_o}$.  Hence $T^*_{\cF}=\pi^*\Omega_E+D_{s_o}$ and 
\begin{eqnarray}
0&=&D_{s_o}\cdot D_{s_o}\nonumber \\
&=&c_1(T^*_{\cF})\cdot c_1(T^*_{\cF})-2c_1(T^*_{\cF})\cdot \pi^*c_1(\Omega_E)+ \pi^*c_1(\Omega_E)\cdot  \pi^*c_1(\Omega_E)\nonumber \\
&=&c_1(T^*_{\cF})\cdot c_1(T^*_{\cF}), \label{eq:zero}
\end{eqnarray}
where we used the fact that $c_1(\Omega_E)=0$.

Now we observe that the foliation is defined on $M$ by the fiber of $\Phi$.  Hence the Chern form $c_1(N^*_{\cF})$ of $N^*_{\cF}$ on $\tM$ is the 
pull-back of a $(1,1)$-form on $C$, which is of dimension $1$.  Hence the pointwise product $c_1(N^*_{\cF})\wedge c_1(N^*_{\cF})$ vanishes
on $M$, which implies that $c_1(N_{\cF})\cdot c_1(N_{\cF})=c_1(N^*_{\cF})\cdot c_1(N^*_{\cF})=0$.  This implies that
\begin{eqnarray*}
0&=&(c_1(T_M)-c_1(T_{\cF}))\cdot (c_1(T_M)-c_1(T_{\cF}))\\
&=&c_1(T_M)\cdot c_1(T_M)-2c_1(T_M)\cdot c_1(T_{\cF})\\
&=&c_1(T_M)\cdot c_1(T_M)-2c_1(N_{\cF})\cdot c_1(T_{\cF})
\end{eqnarray*}
where we have used the identity (\ref{eq:zero}).  This however contradicts the fact that $c_1(T_M)\cdot c_1(T_M)=9$, an odd number.  The contradiction rules out the case of $h^0(M,\Omega_M)=1$ and $h^0(M,N^*_{\cF})=0$.  The rest follows from Lemma 1.

\qed

\ms\ni
{\bf 3.4}  We now proceed to the proof of second main result.

\ms
\ni{\bf Proof of Theorem 2}  

From the classical results of Eells-Sampson and its more recent generalizations, cf. [ES], [L], we know that there exists a $\rho$-equivariant harmonic map
from $\Psi:\tM\rightarrow \tN$, where $\tN\cong B_{\bC}^2$ and $\rho(\Gamma)\subset \Aut(\tN)=PU(2,1)$.
Suppose that $\rank_{\bR}(\Psi)\geqslant 3$, we may apply the result of Siu [Siu1], [Siu2] to conclude that $\Psi$ is  either holomorphic or conjugate-holomorphic as remarked \S1.  In the latter case, we take the complex conjugate of the image so that $\Psi$ becomes holomorphic. In such case, $\Psi$ is a holomorphic mapping of complex rank $2$.

Suppose that $\rank_{\bR}(\Psi)=1$. Then $\Psi$ is a totally geodesic curve $\ell$ in $\tN$ according to a result of Sampson [Sa].  It follows from
our assumption that $\Gamma$ is Zariski dense and the action at $\partial B_{\bC}^2$ does not have fixed points since $\rho$ is reductive, we conclude that it cannot happen that $\rank_{\bR}=1$.  We refer the reader to [CT], \S7
for more elaborations.
Hence it suffices to rule out the case that $\rank_{\bR}(\Psi)=2$.

Let us now study $\Psi$ from a slightly different point of view.  From Lemma 4.5 of [Sim], $\rho$ comes from a complex variation of Hodge structure.
Hence $\rho$ induces an equivariant holomorphic mapping $\Phi:\tM\rightarrow \tS$ such that $\Psi=p\circ\Phi$, where $\tS$ is a Griffiths's period domain
over $N$ and $p:\tS\rightarrow N$ is the projection map.  The only choices of $\tS$ over $B_{\bC}^2=PU(2,1)/P(U(2)\times U(1)$ are 
$PU(2,1)/P(U(2)\times U(1)$ or $PU(2,1)/P(U(1)\times U(1)\times U(1)$.  The first case shows that $\Phi=\Psi$ is holomorphic.  
In the second case, the image $\Phi(\tM)$ lies in a horizontal slice of $PU(2,1)/P(U(1)\times U(1)\times U(1)$ and hence has to lie in
 $PU(2,1)/P(U(2)\times U(1)$ and we conclude that 
 $$\Psi(\tM)\subset PU(2,1)/P(U(2)\times U(1) \subset PU(2,1)/P(U(1)\times U(1)\times U(1).$$
 Hence the fibers of $\Phi$ and hence $\Psi$ are complex curves.  Note that the image has to be a complex curve $C$ in $B_{\bC}^2$.
 By considering normalization if necessary, we may assume that $C$ is a smooth hyperbolic curve with a metric induced from $B_{\bC}^2$.
 In either case, we claim that for
 a mapping $\Psi$ or $\Phi$ as above, all the fibers descend to compact complex curves on $M$.  For simplicity of presentation, we may just regard $\Psi$ as $\Phi$
 in the first case above, and we have a $\rho$-equivariant holomorphic map $\Phi:\tM\rightarrow C$.  Theorem 1 implies that there is a fibration $f:M\rightarrow R$ over a hyperbolic projective algebraic curve $R$ such that $f$ is a fibration.

Now Euler formula for fibrations as in (\ref{eq:euler}) implies that
\begin{equation*}
e(M)=e(R)e(M_s)+\sum_{s_o}n_{s_o}, \label{eq:euler2}
\end{equation*}
As $e(M_s)\geqslant 2$ and $e(R)\geqslant 2$ due to hyperbolicity, we reach a contradiction since $e(M)=3$.  Hence we conclude that $\rank_{\bR}(\Psi)\neq2$.
It follows that $\Phi$ is a holomorphic map of complex rank $2$.

\qed

\ni{\bf Remark} The requirement that $\rho$ is rigid is not really needed, since from Theorem 3 of [Sim], we know that $\rho$ can be deformed to a representation $\rho_o:\Gamma\rightarrow PU(2,1)$ which comes from a complex variation of Hodge structure. 

\bs
\begin{center}
{\bf 4. A rigidity result} 
\end{center} 

\ms\ni
{\bf 4.1} \begin{lemm}
Let $\tM\cong B_{\bC}^2$ and $\tN\cong B_{\bC}^2$.  Let $\Gamma$ be a compact lattice on $PU(2,1)$ and $M=\tM/\Gamma$ with $c_2(M)=3$.
Let $\Phi:\tM\rightarrow \tN$ be a $\Gamma$-equivariant holomorphic map of complex rank $2$.   Then $\Phi$ is a biholomorphism.
\end{lemm}

\ni{\bf Proof} Denote by $C_1(\tN)$ and $C_2(\tN)$ the first and the second Chern forms associated to the Poincar\'e
metric on $\tN$.  Let $\Sigma$ be a fundamental domain of $M$ in $\tM$.  First of all, we claim that 
$$\int_{M}\Phi^*C_1(\tN)\wedge \Phi^*C_1(\tN)=\int_{\Sigma}\Phi^*C_1(\tN)\wedge \Phi^*C_1(\tN)$$
and
$$\int_{M}\Phi^*C_2(\tN)=\int_{\Sigma}\Phi^*C_2(\tN)$$
are positive integers.
The identities above hold since 
the integrants are equivariant under $\Gamma$ and hence descend to $M$.
Clearly the first integral $\int_M \Phi^*C_1(\tN)\wedge \Phi^*C_1(\tN)\in\bZ$, as Chern number of the pull-back line bundle. Now it is well known that on any complex surface $S$ with $p:P(S)\rightarrow S$ the projection map of the projectivized tangent bundle, there is pointwise identification as a differential form
$$\int_{P(M)}\Phi^*C_1(L)^3=-\int_M\Phi^*(C_1^2-C_2).$$
The mapping $\Phi$ induces a mapping between the corresponding projectivized tangent bundles which implies that $\int_{PM}\Phi^*C_1(L)^3\in \bZ$.  Using the identity above
and the fact that  $\int_M\Phi^*C_1^2(\tN)\in \bZ$, we conclude that $\int_M\Phi^*C_2(\tN)\in \bZ$ as well.
Note that the Chern forms of the bundles are positive on $\tN$ in terms of the
Poincar\'e metric there. Hence the integrals are positive as well. The claim is proved.

 In the following, for simplicity of notations, we denote $\Phi^*c_2=\int_M\Phi^*C_2(\tN)$ and
 $\Phi^*c_1^2=\int_M \Phi^*C_1(\tN)\wedge \Phi^*C_1(\tN)$.  It follows that $\Phi^*c_2\in\bZ$.  
As the Chern forms of the Poincar\'e metric on $\tM$ satisfies $C_1^2=3C_2$ pointwise on the form level, we conclude that
$\Phi^*c_1^2=3\Phi^*c_2$ and hence is a multiple of $3$.

Assume that $\Phi$ is not a biholomorphism.  Then the ramification divisor $R$ of $\Phi$ is non-trivial.
We may write $R=\sum_{i=1}^kb_iR_i$, where $R_i$'s are the irreducible components and $b_i$ is the ramification order along $R_i$.
From Riemann-Hurwitz Formula,
\begin{eqnarray}
c_1^2(M)&=&\Phi^*c_1^2(\tN)-2\Phi^*K_{\tN}\cdot R+R\cdot R\label{eq:c1a}\\
&=&\Phi^*c_1^2(\tN)+\Phi^*K_{\tN}\cdot R+K_{\tM}\cdot R\label{eq:c1b}\\
&=&\Phi^*c_1^2(\tN)+\Phi^*K_{\tN}\cdot \sum_{i=1}^kb_iR_i+K_{\tM}\cdot \sum_{i=1}^kb_iR_i.\label{eq:c1c}
\end{eqnarray}
In the above, we have used 
\begin{equation}
K_M=K_{\tM}=\Phi^*K_{\tN}+R \label{eq:RH}.
\end{equation}
From equation (\ref{eq:c1b}), we conclude that $\Phi^*K_{\tN}\cdot R+K_{\tM}\cdot R=6$, where the first term is 
non-negative and the second term is positive.  Note also that $\Phi^*K_{\tN}\cdot R=(-R+K_M)\cdot R$ is even after applying Riemann-Roch to $R$ on $M$.  Hence one of the following holds,

\ms
\ni (a) $\Phi^*K_{\tN}\cdot R=0, K_M\cdot R=6$,\\
(b) $\Phi^*K_{\tN}\cdot R=2 , K_M\cdot R=4$, or\\
(c) $\Phi^*K_{\tN}\cdot R=4 , K_M\cdot R=2$.

For  Case (a), $R$ has to be contracted by $\Phi$ since $K_{\tN}$ is positive from Hodge Index Theorem. This is possible only if $R\cdot R$ is negative from contraction criterion on surface, which would violate equation (\ref{eq:RH}) after taking intersection with $R$. For Case (c), $R\cdot R = K_M \cdot R-\Phi^*K_{\tN}\cdot R=-2$. From Adjunction formula, it follows that the genus of $R$ is $\leqslant1$, which violates the fact that $M$ is hyperbolic. Hence only (b) is possible, with $R\cdot R=2$ as a consequence.

As $-\Phi^*c_1(\tN)\cdot R=2$, it is easy to see that there can at most be two irreducible components in $R$.  First assume that there are two irreducible components $R_1$ and $R_2$ in $R$.  The above constraint leads to $-\Phi^*c_1(\tN)\cdot R_i=1$ and 
$b_i=1$ for $i=1,2$.  For each $i$, $R_i\cdot R_i>0$, for otherwise the Adjunction Formula implies that $R_i$ has genus less
than $2$, contradicting the fact that $M$ is hyperbolic.  It follows that $R_i\cdot R_i=1$ for $i=1,2$ and $R_1\cdot R_2=0$.  But this implies again from Riemann-Roch for $R_1$ that $-\Phi^*c_1(\tN)\cdot R_1=(-c_1(M)-R_2)\cdot R_1$ is even and positive.  Similarly $-\Phi^*c_1(N)\cdot R_2$ is positive and even.  This contradicts the earlier conclusion that
$2=-\Phi^*c_1(\tN)\cdot R=-\Phi^*c_1(\tN)\cdot(R_1+R_2)$.

Hence we may assume that $R=R_1$ has one irreducible component.
It is known that the fundamental domain $\Sigma$ can be taken as a polyhedron and Poincar\'e Polyhedron Theorem holds, with a finite number of faces in the boundary $\pd \Sigma$.  The faces of  $\pd \Sigma$  are identified by actions of elements
$\gamma_1,\cdots,\gamma_k$ in $\Gamma$ to give rise to a compact manifold isometric to $M$. Similarly, the boundary components of $\Phi(\Sigma)$ are identified by induced actions of $\gamma_1,\cdots,\gamma_k$.  By identifying the corresponding boundary components, we get a compact topological space $\Phi(\Sigma)$.   From the earlier discussions, $\Sigma$ is a two-fold branched cover of $\Phi(\Sigma)$.  Since any two fold covering is a Galois covering, $\Phi(\Sigma)$ is obtained from $\Sigma$ by a $\bZ_2$ quotient.

As $R=R_1$ is fixed by an automorphism of $M$ coming from the generator of the $\bZ_2$ action, we conclude that $R$ is totally geodesic. However, for a totally geodesic curve $R_1$ on a complex two ball quotient with possibly self-intersection, we have the formula $K_M\cdot R_1=3(g(\hR_1)-1)$, where $\hR_1$ is the normalization of $R_1$, cf. Lemma 6 of [CKY].  This violates our conclusion in (b) that $K_M\cdot R=4$.  The contradiction implies that $\Phi$ is a biholomorphism and hence Lemma 2.

\qed

\ms
We remark that Domingo Toledo [T] mentioned us that Lemma 4 was known to him with a different
proof.

\ms\ni
{\bf 4.2}  In the following subsections, we would like to use the results above to give details for some unexplained results
in the published version of [Y3].
we recall the following theorem studied in [Y3].

\begin{theo} [Y3]
Let $M=B_{\bC}^2/\Gamma$ be complex two ball quotient with $c_2(M)=3$.  Then $\Gamma$ is arithmetic.
Moreover, $M$ is either a fake projective plane or a Cartwright-Steger surface.  Hence there are altogether $102$ 
such surfaces.
\end{theo} 

We refer to the reader to [Y1], [Y3] for unexplained terminologies in this section.
There are the following main steps in the proof of Theorem.  The first is on the integrality of $\Gamma$, the second is
on Archimedean rigidity of $\Gamma$, and the third is on classification of such surfaces following the work of [PY] and [CS].

We refer the reader to [Y3]  or its updated version in the web for a proof of integrality of $\Gamma$.  Here the image of $\Gamma$
lies in an algebraic group $G$ defined
over a number field $k$ so that with respect to the identity embedding $\sigma_1:k\rightarrow \bC$, the real points of $G$ is
isomorphic to $PU(2,1)$.   

After we know that $\Gamma$ is
integral, to prove that $\Gamma$ is arithmetic, we need to prove Archimedean rigidity statement as given in [Y1] or [Y3], which
amounts to proving that $G^\sigma$, the set of real points of $G$ with respect to another Archimedean embedding of
$k$ in $\bC$, is compact.  The only situation that $G^\sigma$ is non-compact is that $G^\sigma\otimes{\bR}\cong PU(2,1)$
and hence the induced homomorphism $\rho^\sigma:\Gamma\rightarrow G^\sigma$ leads to a harmonic map $\rho^{\sigma}$-equivariant harmonic map $\Phi$ from $\tM$ to $\tN=G^\sigma/K^\sigma=B_{\bC}^2$, where $K^\sigma$ is a maximal compact 
subgroup of $G^\sigma$.  The proof of the second step in [Y3] depends on the fact that $\Phi$ is a holomorphic map, the details
of which is however not provided in [Y3].  The details now follow from Theorem 2 and Lemma 2 obtained in this paper, since the $\rho$ involved
is rigid, following from the fact that the inclusion of $\Gamma$ in $G\otimes \bR$ is locally rigid from Weil's Local Rigidity Theorem, cf. [W].

The third step is a direct consequence of the work of [PY] and [CS], which cover more than just fake projective planes.

\ms\ni
{\bf 4.3}  In this subsection, we elaborate on the statement of [Y3] that there are altogether $102$ smooth complex
two ball quotients with $c_2=3$.  This includes $100$ fake projective planes and
two Cartwright-Steger surfaces, we need to prove that the complex conjugate of a Cartwright-Steger surface is not
biholomorphic to itself.  The corresponding statement was known for fake projective  planes from [KK].  It is stated
in [Y3] that a similar proof works for the Cartwright-Steger surface.  Here we would like to provide the details.
We remark that geometric properties of the Cartwright-Steger surface have been analyzed in [CS] and [CKY].

It is known from the work of Cartwright Steger that the automorphism group of the surface has order $3$. 
From definition, a Cartwright-Steger surface has $h^{1}=1$.  Suppose on the contrary that there is a complex conjugate diffeomorphism $\lambda$ on the surface. It follows that either $\lambda$ or $\lambda^3$ is a conjugate involution, which has a totally real manifold as the fixed point set $F$.  From the first few paragraphs in the proof of Theorem 1 in {\bf 3.3} ,
we know that $h^{1,1}=3$.  Hence it follows from Lefschetz Fixed Point Formula that $F$ contains a component which is one of the followings, sphere $S^2$, real projective plane $P_{\bR}^2$, two-torus $T^2$ or $P_{\bR}\#T^2$ which has a two fold cover that is
$T^2$.   Neither of the first two cases is possible.  In fact, the lift $\tF$ of either set to the universal covering $\tM\cong B_{\bC}^2\subset \bC^2$
has to be compact.  The square of the Euclidean distance function with respect to the origin on $M$, $d_E(0,z)^2=|z_1|^2+|z_2|^2$,
is convex with respect to the Killing metric and hence the restriction of $d_E(0,z)^2$ to $\tF$ is constant by the Maximum Principle. However as $B_{\bC}^2$ is homogeneous, we may identify the origin $0$ of $B_{\bC}^2$ as an arbitrary point on $\tM$ and 
repeat the same argument to conclude that $\tF$ is a point, a contradiction. For the latter two cases, by considering a harmonic map from $T^2$ to $M$ induced by the immersion, the same argument together with Pressman's Theorem (cf. [CE]) leads to a contradiction as well.
Hence there are precisely two Cartwright-Steger surfaces up to biholomorphism. 

\ms\ni
{\bf 4.4}  Here we would also make another correction to the published version of [Y3].  First of all Lemma 1 of [Y3] was 
incorrect and should be discarded.  Originally, Lemma 1 of [Y3] was used in Step two above, but it is no longer needed with the argument of Archimedean rigidity above.  Secondly, \S4 of [Y3], which is
supposed to eliminate the case of $h^1(M)=2$, should be
removed, since Lemma 2 was not correct.  For the proof of the main theorem of [Y3], which is Theorem 3 above, 
it suffices to modify the rest of the argument for Step 1
handle the case of $h^1(M)=2$.  We refer the reader to the last paragraph of {\bf 5.4} of the revised version of [Y3] for details.

\bs
\begin{center}
{\bf 5. Concluding remarks} 
\end{center} 

\bs\ni
{\bf 5.1}  In the study of rigidity problems for complex ball quotients, it is an interesting and important problem to deduce
holomorphicity or anti-holomorphicity of a harmonic map.  The main result of this paper shows that this is possible in
the case that $M$ has the smallest possible Euler number $3$.   It would be interesting to investigate if the restriction to
the smallest Euler number can be removed.  It is not known even if we allow $c_2(M)$ to be $3n$ with $n\neq 1$ small
as mentioned in Proposition 1.  The case for cofinite complex two ball quotients has not been studied as well.  In such case,
the possible Euler numbers take the values of all possible integers as given in Proposition 2.  At this point, Question 2 and its
analogue for cofinite lattice  of $PU(n,1)$, are essentially completely open.

\ms\ni
{\bf 5.2}  Comparing to known techniques in the study of rigidity problems, the main contributions in this paper is
perhaps the study of the foliation in terms of constraints given by Chern numbers.  In general we lack a good criterion
to force algebraicity of leaves of foliation.  Any non-trivial results in this direction may give significant implications for
different problems in geometry.

\ms\ni
{\bf 5.3} For Question 1, to the knowledge of the author, not much is known in general.  Known techniques in this direction include the work of [CT] and [Sim].
We may confine our discussions to the situations that $M$ and $N$ are compact.
It would be desirable to explore the borderline cases with $\rank_{\bR}=2p(N)$ first.

We remark that if $M$ is a locally Hermitian symmetric spaces of rank at least $2$
and $N$ is a manifold of non-positive Riemannian sectional curvature, such a mapping $\Phi$ has to be a totally geodesic isometry, following 
geometric superrigidity argument as in [MSY].  Note that a holomorphic map between K\"ahler manifolds is necessarily harmonic with
respect to the K\"ahler metrics.

\ms\ni
{\bf 5.4} Related to the earlier problems, we may also consider the problem of
distinguishing a holomorphic isometry from a totally geodesic isometry
as a  mapping between two bounded Hermitian symmetric domains.
For this active area of research, an interesting reader may consult [Mo3] for a survey and references of recent results.

\bigskip

\noindent{\bf References}

\bs
\ni [AMRT] Ash A., Mumford D., Rapoport M., Tai Y-S., Smooth Compactifications of Locally Symmetric Varieties,
Math. Sci. Press 1975.

\ms
\ni [BHPV] Barth, W. P., Hulek, K., Peters, C. A. M., Van de
Ven, A., Compact complex surfaces. Second edition. Ergebnisse der
Mathematik und ihrer Grenzgebiete. 3. Folge. A Series of Modern
Surveys in Mathematics 4. Springer-Verlag, Berlin, 2004.

\ms
\ni [BB] Baum, P., Bott. R., On the zeroes of meromorphic vector fields,  Essays on Topology and Related Topics, Springer-Verlag,
New York, 1970, p. 29-47.

\ms
\ni [Br1] Brunella, M., Feuilletages holomorphes sur les surfaces complexes compactes. Ann. Sci. cole Norm. Sup. 30 (1997), 569-594.

\ms
\ni [Br2] Brunella, M., Birational geometry of foliations. IMPA Monographs, 1. Springer, Cham, 2015.

\ms
\ni [CKY] Cartwright, D., Koziarz, V. and Yeung, S.-K., On the Cartwright-Steger surface, arXiv:1412.4137, to appear in J. Alg. Geom.

\ms
\ni [CS] Cartwright, D., Steger, T.,  Enumeration of the $50$ fake projective planes, C. R. Acad. Sci. Paris, Ser. 1,
348 (2010), 11-13,  see also \\
\verb'http://www.maths.usyd.edu.au/u/donaldc/fakeprojectiveplanes/'

\ms
\ni [CT]  Carlson, J. A., Toledo, D., Harmonic mappings of K\"ahler manifolds to locally symmetric spaces. Inst. Hautes \'Etudes Sci. Publ. Math. No. 69 (1989), 173-201.

\ms
\ni [CE]   Cheeger, J.,  Ebin, D. G., Comparison theorems in Riemannian geometry. North-Holland Mathematical Library, Vol. 9. North-Holland Publishing Co., Amsterdam-Oxford; American Elsevier Publishing Co., Inc., New York, 1975.

\ms
\ni [ES] Eells, J. and Sampson, J., 
Harmonic mappings of Riemannian manifolds, Amer. J. Math. 86 (1964), 
109-160.

\ms
\ni [H] Hirzebruch, F., Chern numbers of algebraic surfaces: an example, Math. Ann. 266 (1984), 351-356.

\ms
\ni  [JY] Jost, J., Yau, S.-T., Harmonic maps and group representations. Differential geometry, 241-259, Pitman Monogr. Surveys Pure Appl. Math., 52, Longman Sci. Tech., Harlow, 1991.

\ms
\ni[KK] Kharlamov,\:V., Kulikov,\:V., On real structres on rigid surfaces.  Izv.\,Math. {66} (2002), 133-150.

\ms
\ni [L] Labourie, F., 
Existence d'applications  harmoniques tordues \'a 
valeurs dans les vari\'et\'es \'a courbure
n\'egative,  Proc. Amer. Math. Soc. 111 (1991), 877-882.

\ms
\ni [Mo1] Mok, N.,  The holomorphic or antiholomorphic character of harmonic maps into irreducible compact quotients of polydiscs. Math. Ann. 272 (1985), 197-216.

\ms
\ni [Mo2] Mok, N., Projective-algebraicity of minimal compactifications of complex-hyperbolic space forms of finite volume, in: Perspectives in analysis, geometry, and topology. On the occasion of the 60th birthday of Oleg Viro, Progr. Math. 296 (2012), Birkhuser-Verlag, Basel (2008), 331-354.
 
\ms
\ni [Mo3] Mok, N.,  Geometry of holomorphic isometries and related maps between bounded domains. Geometry and analysis. No. 2, 225-270, Adv. Lect. Math. (ALM), 18, Int. Press, Somerville, MA, 2011.

\ms
\ni Mok, N., Siu, Y.-T., Yeung, S.-K., Geometric superrigidity, Invent. Math., 113 (1993), pp. 57-83.

\ms
\ni [Mu] Mumford, D., Hirzebruch's Proportionality Theorem in the
non-compact case, Inv. Math. 42 (1977), 239-272.

\ms
\ni [P] Parker, J. R., On the volumes of cusped, complex hyperbolic manifolds and orbifolds,
Duke Math. J. 94 (1998), 433-464.

\ms
\ni [PY] Prasad, G., and Yeung, S.-K.,  Fake projective planes. Inv.\,Math. 168 (2007), 321-370; Addendum, ibid 182 (2010), 213-227.

\ms
\ni [Sa] Sampson, J. H., Some properties and applications of harmonic mappings, Ann. Sci.\'Ecole Norm. Sup., 11 (1978), 211-228.

\ms
\ni [Sim] Simpson, C., Higgs bundles and local systems. Inst. Hautes \'Etudes Sci. Publ. Math. No. 75 (1992), 5-95.

\ms
\ni [Siu1] Siu, Y.-T,  The complex-analyticity of harmonic maps and the strong
rigidity of compact K\"ahler manifolds, Ann. of Math. 112 (1980), 73-111.

\ms
\ni [Siu12 Siu, Y.-T,  
Complex-analyticity of harmonic maps, vanishing and Lefschetz theorems. J. Differential Geom. 17 (1982), 55-138.

\ms
\ni [Siu3] Siu, Y.-T., Strong rigidity for K\"ahler manifolds and the construction of bounded holomorphic functions, in discrete groups in
geometry and analysis, r. Howe (editor), Birkhauser, 1987, 124-151.


\ms
\ni [T] Toledo, D., private communication.

\ms
\ni [W] Weil, A., Discrete subgroups of Lie groups II, Ann. of Math., 75 (1962), pp. 578-602.

\ms
\ni [Y1] Yeung, S.-K., Integrality and arithmeticity of co-compact lattices corresponding to certain complex two
ball quotients of Picard number one, Asian J. Math. 8 (2004), 104-130; Erratum, Asian J. Math. 13
(2009), 283-286.

\ms
\ni [Y2]  Yeung, S.-K., Classification and construction of fake projective planes, Handbook of geometric analysis, No. 2, 391-431, Adv. Lect. Math. (ALM), 13, Int. Press, Somerville, MA, 2010.

\ms
\ni [Y3] Yeung, S.-K., Classification of surfaces of general type with Euler number $3$, J. Reine Angew. Math. 697 (2014), 1-14, corrected version,\\
\verb'http://www.math.purdue.edu/\%7Eyeung/papers/integralpoints-16b.pdf/'

\end{document}